\documentclass[reqno,fleqn,onecolumn]{amsproc}%
\usepackage{amsfonts}
\usepackage{amsmath}
\usepackage{amssymb}
\usepackage{graphicx}
\usepackage[none]{hyphenat}%
\theoremstyle{plain}

\numberwithin{equation}{section}
\begin{document}
\title{}

\begin{center}
{\LARGE \textbf{Delta Hedging without the Black-Scholes Formula}}

\bigskip

Yukio Hirashita

\bigskip
\end{center}

\noindent\textbf{Abstract: }We introduce a new method of delta hedging. In
many cases, this method results in a lower cost than the Black-Scholes method.
To calculate the cost of hedging, we develop a Mathematica program.

\bigskip



\bigskip

\noindent{\large \textbf{1. Introduction}}

\noindent We introduce a pricing method of distribution (see Hirashita [4])
under the condition that the risk-free interest rate is equal to the growth
rate. This pricing method yields a new method of delta hedging. In this paper,
we re-calculate several examples given in Hull [5] and Luenberger [7], and
show that in many cases, the new method of hedging results a lower cost than
the Black-Scholes method. In order to calculate the cost of hedging under this
new method, we utilize a Mathematica program that include the two-dimensional
Newton-Raphson method. It should be noted that Mathematica is a programming language.

\bigskip

\noindent{\large \textbf{2. Well-known results: Delta and gamma with the
Black-Scholes formula}}

\noindent We assume that the stock price $Y=Se^{x+rT}$ is lognormally
distributed with volatility $\sigma\sqrt{T}$, where $S$ is the stock price
factor, $r$ is the continuously compounded interest rate, $K$ is the exercise
price of the call option, and $T$ is the exercise period. Then, the European
call option is given by
\begin{align*}
a(x)  &  :=\max(Se^{x+rT}-K,\text{ }0)=\left\{
\begin{array}
[c]{cc}%
Se^{x+rT}-K, & \text{if }x\geq\log\frac{K}{S}-rT,\\
0, & \text{others,}%
\end{array}
\right. \\
dF(x)  &  :=p(x)dx:=\frac{1}{\sqrt{2\pi T}\sigma}e^{-\frac{(x+\sigma
^{2}T/2)^{2}}{2\sigma^{2}T}}dx,
\end{align*}
and $x\in I:=(-\infty,+\infty)$. The expectation $E$ of this option is
calculated as
\begin{align*}
E  &  =\int_{I}a(x)dF(x)=\frac{1}{\sqrt{2\pi T}\sigma}\int_{\log\frac{K}%
{S}-rT}^{\infty}(Se^{x+rT}-K)e^{-\frac{(x+\sigma^{2}T/2)^{2}}{2\sigma^{2}T}%
}dx\\
&  =Se^{rT}N\left(  \frac{\log\frac{S}{K}+(r+\frac{\sigma^{2}}{2})T}%
{\sigma\sqrt{T}}\right)  -KN\left(  \frac{\log\frac{S}{K}+(r-\frac{\sigma^{2}%
}{2})T}{\sigma\sqrt{T}}\right)  ,
\end{align*}
where $N(x)=\int_{-\infty}^{x}e^{-x^{2}/2}/\sqrt{2\pi}\,dx$ is the cumulative
standard normal distribution function. We set%
\[
d_{1}=\frac{\log\frac{S}{K}+(r+\frac{\sigma^{2}}{2})T}{\sigma\sqrt{T}}\text{
\ and \ }d_{2}=\frac{\log\frac{S}{K}+(r-\frac{\sigma^{2}}{2})T}{\sigma\sqrt
{T}};
\]
then, the equation $E/\underline{u}=e^{rT}$ yields the price%
\begin{equation}
\underline{u}=SN\left(  d_{1}\right)  -Ke^{-rT}N\left(  d_{2}\right)  ,
\tag{2.1}%
\end{equation}
which is the Black-Scholes formula for a European call option (see Hull [5]).
Delta and gamma are given as follows:%

\begin{equation}
\text{delta}:\underline{\Delta}:=\frac{\partial\underline{u}}{\partial
S}=N\left(  d_{1}\right)  . \tag{2.2}%
\end{equation}

\begin{equation}
\text{gamma}:\underline{\Gamma}:=\frac{\partial^{2}\underline{u}}{\partial
S^{2}}=e^{-d_{1}^{2}/2}/(S\sigma\sqrt{2\pi T}). \tag{2.3}%
\end{equation}

\bigskip

\noindent{\large \textbf{3. New results: Delta and gamma with the simultaneous
equations}}

\noindent It should be noted that $\inf_{x\in I}$ $a(x)=0$ and $\int
_{I}a(x)dF(x)<\infty.$ As $\int_{a(x)=0}dF(x)>$ $0,$ the price $u$ and the
optimal proportion of investment $t_{u}$ are determined by the simultaneous
equations%
\begin{equation}
\left\{
\begin{array}
[c]{c}%
\exp(\int_{I}\log(\frac{a(x)t_{u}}{u}-t_{u}+1)dF(x))=e^{rT},\\
\int_{I}\frac{a(x)-u}{a(x)t_{u}-ut_{u}+u}dF(x)=0
\end{array}
\right.  \tag{3.1}%
\end{equation}
(see Corollary 5.1 and Section 6 in Hirashita [4]). Set $\beta(x):=a(x)t_{u}%
-ut_{u}+u$; then, we obtain%
\[
\left\{
\begin{array}
[c]{c}%
\int_{I}\log\beta(x)dF(x)=rT+\log u,\\
\int_{I}\frac{1}{\beta(x)}dF(x)=\frac{1}{u}.
\end{array}
\right.
\]
As $a(x)$, $t_{u},$ and $u$ are functions with respect to $S$, from
$\partial\left(  \int_{I}\log\beta(x)dF(x)-\log u\right)  /\partial S$ $=0,$
we have
\[
t_{u}\int_{I}\frac{\frac{\partial a}{\partial S}(x)}{\beta(x)}dF(x)+\frac
{\partial t_{u}}{\partial S}\int_{I}\frac{a(x)-u}{\beta(x)}dF(x)+(1-t_{u}%
)\frac{\partial u}{\partial S}\int_{I}\frac{1}{\beta(x)}dF(x)=\frac{1}{u}%
\frac{\partial u}{\partial S},
\]
this implies that%
\begin{equation}
\text{delta: }\Delta:=\frac{\partial u}{\partial S}=u\int_{I}\frac
{\frac{\partial a}{\partial S}(x)}{\beta(x)}dF(x)=uWe^{rT}, \tag{3.2}%
\end{equation}
where%
\[
W:=\int_{\log\frac{K}{S}-rT}^{\infty}\frac{e^{x}}{St_{u}e^{x+rT}-Kt_{u}%
-ut_{u}+u}p(x)dx.
\]

\noindent Using the well-known formula%
\[
\frac{\partial}{\partial S}\int_{\varphi(S)}^{\psi(S)}F(S,x)dx=F(S,\psi
(S))\psi^{\prime}(S)-F(S,\varphi(S))\varphi^{\prime}(S)+\int_{\varphi
(S)}^{\psi(S)}\frac{\partial F}{\partial S}(S,x)dx,
\]
we have
\begin{align*}
\frac{\partial W}{\partial S}  &  =\frac{K}{u(1-t_{u})S^{2}e^{rT}}p(\log
\frac{K}{S}-rT)\\
&  -\int_{\log\frac{K}{S}-rT}^{\infty}\frac{t_{u}e^{x+rT}+\left(
Se^{x+rT}-K-u\right)  \frac{\partial t_{u}}{\partial S}+(1-t_{u})uWe^{rT}%
}{\left(  St_{u}e^{x+rT}-Kt_{u}-ut_{u}+u\right)  ^{2}}e^{x}p(x)dx.
\end{align*}
From $\partial\left(  \int_{I}1/\beta(x)dF(x)-1/u\right)  /\partial S$ $=0$,
we obtain
\begin{align*}
&  \frac{\partial t_{u}}{\partial S}\\
&  =\frac{(1-t_{u})uWe^{rT}\int_{I}\frac{1}{\beta^{2}(x)}dF(x)+t_{u}e^{rT}%
\int_{\log\frac{K}{S}-rT}^{\infty}\frac{e^{x}}{\left(  \left(  Se^{x+rT}%
-K\right)  t_{u}-ut_{u}+u\right)  ^{2}}dF(x)-\frac{We^{rT}}{u}}{-\int_{I}%
\frac{a(x)-u}{\beta^{2}(x)}dF(x)}.
\end{align*}
Therefore, we can calculate%
\begin{equation}
\text{gamma: }\Gamma:=\frac{\partial^{2}u}{\partial S^{2}}=\frac{\partial
}{\partial S}(uWe^{rT})=uW^{2}e^{2rT}+ue^{rT}\frac{\partial W}{\partial S}.
\tag{3.3}%
\end{equation}

\bigskip

\noindent{\large \textbf{4. Comparison between two delta hedging methods}}

\noindent When $S=49$, $K=50$, $\sigma=0.2$, $r=0.05$, and $T=20/52$ (see
Section 14.1 in Hull [5]), we have%
\[
a(x)=\max(49e^{1/52}e^{x}-50,\text{ }0),\text{ \ \ }dF(x)=\sqrt{\frac{65}%
{2\pi}}e^{-\frac{65(x+\frac{1}{130})^{2}}{2}}dx,
\]
and the price $u\doteq1.774$ with respect to Equation 3.1. At this price, if
investors continue to invest $t_{u}$ $\doteq0.115$\ of their current capital,
they can maximize the limit expectation of the growth rate to $e^{rT}%
=e^{1/52}\doteq1.019.$ Moreover, we have $\Delta\doteq0.448$ and $\Gamma
\doteq0.0668.$

Meanwhile, the Black-Scholes formula yields $\underline{u}\doteq2.401$,
$\underline{\Delta}\doteq0.522,$ and $\underline{\Gamma}\doteq0.0655$.

\bigskip

In the following examples, we consider the cost of hedging for $100,000$ stocks.

\bigskip

\textbf{Example 4.1}. Table 14.2 in Hull [5] presents a simulation of the
delta hedging of a sequence of weekly stock prices $S_{1}:=\{49.00,$ $48.12,$
$47.37,$ $50.25,$ $51.75,$ $53.12,$ $53.00,$ $51.87,$ $51.38,$ $53.00,$
$49.88,$ $48.50,$ $49.88,$ $50.37,$ $52.13,$ $51.88,$ $52.87,$ $54.87,$
$54.62,$ $55.87,$ $57.25\}$ with $K=50$, $\sigma=0.2,$ and $r=0.05$. For
$S_{1}$, the cost of hedging is $\$287,500$ according to Equation 3.2, which
is $9.2\%$ higher than the cost of hedging using the Black-Scholes method
($\$263,300$).

On the other hand, when $K=65$, the cost of hedging is $\$2,600$ according to
Equation 3.2, which is $46.9\%$ lower than the cost of hedging using the
Black-Scholes method ($\$4,900$).

In addition, when $K=35$, the difference between these costs is within $0.1\%$.

\bigskip

\textbf{Example 4.2.} Table 14.3 in Hull [5] shows a simulation of the delta
hedging of a sequence of weekly stock prices $S_{2}:=\{49.00$, $49.75,$
$52.00$, $50.00$, $48.38$, $48.25$, $48.75$, $49.63$, $48.25$, $48.25$,
$51.12$, $51.50$, $49.88$, $49.88$, $48.75$, $47.50$, $48.00$, $46.25$,
$48.13$, $46.63$, $48.12\}$ with $K=50$, $\sigma=0.2,$ and $r=0.05$. For
$S_{2}$, the cost of hedging is $\$247,900$ according to Equation 3.2, which
is $3.4\%$ lower than the cost of hedging using the Black-Scholes method
($\$256,600$).

On the other hand, when $K=65$, the cost of hedging is $\$3,100$ according to
Equation 3.2, which is $48.3\%$ lower than the cost of hedging using the
Black-Scholes method ($\$6,000$).

In addition, when $K=35$, the difference between these costs is within $0.1\%$.

\bigskip

\textbf{Example 4.3}. We select Table 13.1 in Luenberger [7] of weekly stock
prices\ such that $S_{3}:=\{35.50,$ $34.63,$ $33.75,$ $34.75,$ $33.75,$
$33.00,$ $33.88,$ $34.50,$ $33.75,$ $34.75,$ $34.38,$ $35.13,$ $36.00,$
$37.00,$ $36.88,$ $38.75,$ $37.88,$ $38.00,$ $38.63,$ $38.50,$ $37.50\}$ with
$K=35$, $\sigma=0.18,$ and $r=0.1$. For $S_{3}$, the cost of hedging is
$\$274,900$ according to Equation 3.2, which is $1.3\%$ higher than the cost
of hedging using the Black-Scholes method ($\$271,300)$.

On the other hand, when $K=45$, the cost of hedging is $\$1,300$ according to
Equation 3.2, which is $82.7\%$ lower than the cost of hedging using the
Black-Scholes method ($\$7,500$).

In addition, when $K=25$, the difference between these costs is within $0.1\%$.

\bigskip

\noindent{\large \textbf{5. Program}}

\noindent In order to confirm that the cost of hedging is indeed $\$287,500$
in Example 4.1, it is necessary to run the following Mathematica program that
includes the Newton-Raphson method. In order to confirm the other costs of
hedging with respect to Equation 3.2, substitute the \texttt{Stock},
\texttt{K}, \texttt{sigma}, and \texttt{r} data with those given below.

\bigskip
\begin{verbatim}
     Stock = {49, 48.12, 47.37, 50.25, 51.75, 53.12, 53, 51.87, 51.38,
     53, 49.88, 48.50, 49.88, 50.37, 52.13, 51.88, 52.87, 54.87, 54.62,
     55.87, 57.25};
     K = 50; sigma = 0.2; r = 0.05;
     m = 20; shares = 100000;
     a[x_] := If[x < Log[K/S] - r*T, 0, S*Exp[r*T]*Exp[x] - K];
     p[x_] := Exp[-(x + sigma^2*T/2)^2/(2*sigma^2*T)]/Sqrt[2Pi*T*sigma^2];
     wr = 50; mr = 10; sd = 1000; M = 10^12;
     Price := Module[{u, t},
     u = NIntegrate[a[x]p[x], {x, -Infinity, Infinity}, WorkingPrecision
     -> wr, MaxRecursion -> mr, SingularityDepth -> sd]/2; t=0.5;
     Do[{If[u < 1/M, Break[]];
       f = Exp[NIntegrate[Log[(a[x]t/u - t + 1)]*p[x], {x, -Infinity,
       Infinity}, WorkingPrecision -> wr, MaxRecursion -> mr,
       SingularityDepth -> sd]] - Exp[r*T];
       g = NIntegrate[(a[x] - u)/(a[x]t - u*t + u)*p[x], {x, -Infinity,
       Infinity},  WorkingPrecision -> wr, MaxRecursion -> mr,
       SingularityDepth -> sd];
       fu = (f + Exp[r*T])*NIntegrate[-a[x]t/(a[x]t - u*t + u)/u*p[x],
       {x, -Infinity, Infinity}, WorkingPrecision -> wr, MaxRecursion ->
       mr, SingularityDepth -> sd];
       ft = 0;
       gu = -NIntegrate[a[x]/(a[x]t - u*t + u)^2*p[x], {x, -Infinity,
       Infinity}, WorkingPrecision -> wr, MaxRecursion -> mr,
       SingularityDepth -> sd];
       gt = -NIntegrate[(a[x] - u)^2/(a[x]t - u*t + u)^2*p[x], {x,
       -Infinity, Infinity}, WorkingPrecision -> wr, MaxRecursion -> mr,
       SingularityDepth -> sd];
       ans = Solve[Re[fu]*a + Re[ft]*b == -Re[f] && Re[gu]*a + Re[gt]*b
       == -Re[g], {a, b}]; If[ans=={},Break[]];a0 = a /. ans[[1]];
       b0 = b /. ans[[1]]; If[Abs[a0] + Abs[b0] < 1/M, Break[]]; u = u + a0;
       t = t + b0; If[u < 0, u = (u - a0)/2]; If[t < 0, t = (t - b0)/2];
       If[t >= 1, t = ((t - b0)+1)/2];
        }, {j2, 1, 100}]; t0 = t;
     u];
     Print["{Week, Stock price, Delta, Shares, Cost of shares,
     Cumulative cost, Interest}"];
     s = 0; cost = 0; interest = 0;
     Do[S = Stock[[j + 1]]; T = (m - j)/52; If[T > 0, u = Price; t = t0,
       u = If[S > K, S - K, 0]; t =.];
       If[T > 0, W = NIntegrate[Exp[x]/((S*Exp[r* T]*Exp[x] - K)t -
       u*t + u)*p[x], {x, Log[K/S] - r*T, Infinity}, WorkingPrecision ->
       wr,  MaxRecursion -> mr, SingularityDepth -> sd], If[ S > K,
       W = 1/(S - K), W = 0]];
       delta = u*Exp[r*T]*W; delta = Round[1000*delta]/1000.;
       s2 = shares*delta - s; s = shares*delta;
       cost = cost + Round[s2*S/100]/10 + interest;
       interest = Round[cost*r*1/52*10]/10.;
       Print[{j, S, u, delta, s2, Round[s2*S/100]/10., Round[cost*1000],
       interest}];,
     {j, 0, m}];
     cost = cost - delta*shares*Min[K, Stock[[m + 1]]]/1000 + (1 - delta)*
     shares*Max[Stock[[m + 1]] - K, 0]/1000;
     Print["K = ", K, ", Cost of Hedging =", Round[cost*1000]];
\end{verbatim}

\bigskip

\bigskip

\noindent\textbf{References}

\bigskip

\noindent\lbrack1] D. Bernoulli, Exposition of a new theory of the measurement of\ risk,

Econometrica 22 (1954), 23--36.

\noindent\lbrack2] F. Black and M. Scholes, The pricing of options and
corporate liabilities,

Journal of Political Economy 81 (1973), 637-54.

\noindent\lbrack3] W. Feller, An introduction to probability theory and \ its
application, John

Wiley and Sons, New York, 1957.

\noindent\lbrack4] Y. Hirashita, Game pricing and double sequence of random variables,

Preprint (2007), arXiv:math.OC/0703076.

\noindent\lbrack5] J. Hull, Options, futures, and other derivatives, Prentice
Hall, New Jersey,

2003.

\noindent\lbrack6] J. L. Kelly, A new interpretation of information
rate,\ Bell system Technical

Journal 35 (1956), 917--926.

\noindent\lbrack7] D. G. Luenberger, Investment science, Oxford University
Press, Oxford, 1998.

\bigskip

Chukyo university, Japan

E-mail address: yukioh@cnc.chukyo-u.ac.jp

\end{document}